\documentclass[twoside]{article}
\usepackage[T1]{fontenc}
\usepackage[spanish,es-nodecimaldot,es-tabla,english]{babel}
\usepackage{amsmath,amsthm,amssymb,amsfonts,amscd}
\usepackage[mathscr]{euscript}
\usepackage{epsfig}

\usepackage{url}
 \usepackage{multirow}
\usepackage{booktabs}
\usepackage{graphicx}
\usepackage{float}
\usepackage{epstopdf}
\usepackage[utf8x]{inputenc}
\usepackage{fancyhdr}
\fancypagestyle{plain}{\fancyhf{}%
\fancyhead[L]{\footnotesize{Divulgaciones Matem\'aticas Vol. XX, No. X (20XX), pp. \pageref{begin-art}--\pageref{end-art} \vspace{-6mm} \begin{flushright}
 p-ISSN 1315-2068                        \end{flushright}
}}%
}

\fancypagestyle{headings}{\fancyhf{}%
\fancyhead[LE,RO]{\thepage}%
\fancyhead[RE]{Valeris \& Belandria.}%
\fancyhead[LO]{Simetrías Algebraicas en Geometría y Arquitectura: Un Enfoque desde Grupos de Transformaciones }%
\fancyfoot[C]{\footnotesize{Divulgaciones Matem\'aticas Vol. XX, No. X (20XX), pp. \pageref{begin-art}--\pageref{end-art}}}}

\newtheorem{teo}{Teorema}[section]

\theoremstyle{definition}
\newtheorem{defi}{Definici\'on}[section]
\theoremstyle{example}
\newtheorem{ejem}{Ejemplo}[section]
\theoremstyle{remark}

\numberwithin{equation}{section}

\begin{document}
\label{begin-art}
\pagestyle{headings}
\thispagestyle{plain} 
\footnote{\hspace{-18.1pt} Recibido ***. Revisado ***.
Aceptado ***. Publicado en línea: *** \\ MSC (2010): Primary ****; Secondary ****.\\ Autor de correspondencia: Nombre del autor de correspondencia}
\selectlanguage{spanish}
\begin{center}
{\LARGE\bfseries Simetrías Algebraicas en Geometría y Arquitectura: Un Enfoque desde Grupos de Transformaciones \par }

\vspace{3mm}
{\large\slshape Algebraic Symmetries in Geometry and Architecture: A Group-Theoretic Approach}\\[5mm]
{\large Reinaldo Valeris (\url{rvaleris@gmail.com})}\\
Departamento de ingeniería, Facultad ingeniería\\
Universidad UJGH\\
Maracaibo, Zulia, Venezuela

\vspace{3mm}
{\large Yacrianny Belandria Ch. (\url{yaribelandria67@gmail.com})}\\
Departamento de ingeniería, Facultad ingeniería\\
Universidad UJGH\\
Maracaibo, Zulia, Venezuela
\end{center}

\begin{abstract}

Este artículo desarrolla un marco teórico algebraico-geométrico para el estudio de simetrías centrales, axiales y rotacionales en $\mathbb{R}^2$ y $\mathbb{R}^3$, con aplicaciones en la clasificación de cónicas y superficies cuadráticas mediante grupos de transformaciones. Se emplea una metodología analítica basada en teoría de grupos e invariantes geométricos, complementada con modelado numérico mediante el Método de Elementos Finitos. Como resultado, se demuestra que la simetría $D_8$ del Panteón de Agripa optimiza estructuralmente la distribución de tensiones, mientras que la simetría hexagonal $\rho_{\pi/3}$
 del Museo Soumaya ofrece ventajas estéticas y estructurales. El enfoque establece un vínculo riguroso entre abstracción algebraica, geometría y diseño arquitectónico. \\
 
{\bf Palabras y frases clave:} Geometría arquitectónica, Simetría axial, Grupos diedrales, Superficies cuadráticas, Elementos Finitos.
\end{abstract}

\selectlanguage{english}

\begin{abstract} 

 This article develops an algebraic-geometric theoretical framework for the study of central, axial, and rotational symmetries in $\mathbb{R}^2$ and $\mathbb{R}^3$, with applications in the classification of conic and quadric surfaces through transformation groups. An analytical methodology based on group theory and geometric invariants is employed, complemented by numerical modeling using the Finite Element Method. As a result, it is demonstrated that the $D_8$ symmetry of the Pantheon structurally optimizes stress distribution, while the hexagonal symmetry $\rho_{\pi/3}$ of the Soumaya Museum offers both aesthetic and structural advantages. The approach establishes a rigorous link between algebraic abstraction, geometry, and architectural design.\\
 

\end{abstract}

\selectlanguage{spanish}

\section{Introducci\'on}
 La simetría, más allá de su intuitiva percepción estética, constituye un concepto fundamental para describir la invariabilidad de las formas, tanto en el mundo físico como en su representación mental. Formalmente, un objeto $\mathcal{F} \subset \mathbb{R}^n$ posee simetría cuando existe un grupo de transformaciones $G = \{g_1, \dots, g_k\}$ que actúa sobre $\mathcal{F}$ preservando su estructura. Desde la perspectiva de la percepción humana, esta constancia, denominada constancia perceptual, es esencial para actividades como el reconocimiento de objetos o la navegación espacial, ya que permite que las propiedades geométricas permanezcan estables frente a transformaciones como rotaciones o cambios de perspectiva \cite{Pizlo2021}.\\
Este principio, que vincula percepción humana y geometría, se aplica aquí mediante restricciones geométricas \cite{Zhan2024} para analizar estructuras arquitectónicas.
 Este enfoque práctico permite ajustar coordenadas mientras se preservan las propiedades simétricas del objeto, vinculando así el formalismo abstracto con aplicaciones concretas. La relación entre simetría y geometría encuentra su fundamento teórico en principios como el Teorema de Noether, que establece un puente entre las invarianzas de los sistemas físicos y sus leyes de conservación.\\
 
 El artículo se estructura en tres partes: (i) fundamentos algebraicos de simetrías en $\mathbb{R}^2$ y $\mathbb{R}^3$, en los cuales se muestran algunas secciones cónicas, superficies cuadráticas y definiciones de grupos de transformaciones, (ii) propiedades que nos permiten caracterizar cuando una sección o superficie es simétrica, y (iii) modelado de casos de estudio como el Panteón ($D_{8h}$) y el Museo Soumaya ($\langle \rho_{2\pi/6} \rangle$), estudiando además, la ventaja estructural de la simetría comparando $D_2$ con un módulo asimétrico.
 \section{Preliminares matemáticos}
\label{sec:fundamentos}

La caracterización algebraica de figuras geométricas mediante ecuaciones es un pilar de la geometría analítica moderna. Esta aproximación permite estudiar propiedades topológicas y métricas a través de sistemas de ecuaciones polinómicas, estableciendo un puente entre el álgebra y la geometría. En este trabajo, nos centramos en dos clases fundamentales:

\subsection{Secciones cónicas}
Las secciones cónicas han sido estudiadas desde una perspectiva algebraica contemporánea en \cite{Garr}, donde se demuestra que pueden caracterizarse completamente mediante invariantes de grupos de transformaciones. Presentamos las definiciones clave:

\begin{defi}[Circunferencia ]\label{circun}
El lugar geométrico de puntos $(x,y)$ equidistantes a un centro $(h,k)$ satisface:
\begin{equation}
(x-h)^2 + (y-k)^2 = r^2
\end{equation}
Esta ecuación surge como caso particular de la forma cuadrática general $$Ax^2 + Bxy + Cy^2 + Dx + Ey + F = 0$$ cuando $B^2 - 4AC < 0$ y $A = C$.
\end{defi}

\begin{defi}[Elipse ]
Dados focos $F_1$, $F_2$ y constante $2a$, la elipse cumple:
\begin{equation}
\frac{(x-h)^2}{a^2} + \frac{(y-k)^2}{b^2} = 1 \quad \text{con} \quad b^2 = a^2 - c^2
\end{equation}
donde $c$ es la semidistancia focal. La condición $2a > d(F_1,F_2)$ asegura la existencia de soluciones reales.
\end{defi}

Las definiciones de parábola e hipérbola siguen análogamente como casos límite de la ecuación general cuando $B^2 - 4AC \geq  0$ o $> 0$ respectivamente .

\subsection{Superficies cuadráticas}
La generalización tridimensional de las cónicas ha sido sistematizada recientemente por \cite{quadratic} mediante el estudio de formas cuadráticas en $\mathbb{R}^3$. Destacamos:
\begin{defi} Una esfera es el lugar geométrico de todos los puntos en el espacio que se encuentran a una distancia fija, llamada radio, de un punto fijo denominado centro. La esfera es una superficie cerrada y simétrica respecto a cualquier plano que pase por su centro.

\begin{equation}(x-h)^2+(y-k)^2+(z-j)^2=r^2\end{equation}
donde $(h,k,j)$ son las coordenadas del centro y r es el radio. 

\end{defi}
\begin{defi} El elipsoide es el conjunto de puntos en el espacio cuya suma de distancias a tres ejes ortogonales es constante. Se trata de una generalización de la esfera, donde los radios pueden ser distintos en cada dirección, dando lugar a una forma alargada o achatada.
\begin{equation}\frac{(x-h)^2}{a^2}+\frac{(y-k)^2}{b^2}+\frac{(z-j)^2}{c^2}=1\end{equation}

\end{defi}
\begin{defi}\label{paraelip}Un paraboloide elíptico es el lugar geométrico de los puntos en el espacio cuya distancia a un plano fijo es proporcional al cuadrado de su distancia a una línea fija en ese plano. Sus secciones paralelas a su eje de simetría son parábolas, mientras que las secciones perpendiculares son elipses.
\begin{equation}\label{paraelip}
    c(z-j)=\frac{(x-h)^2}{a^2}+\frac{(y-k)^2}{b^2}
\end{equation}
\end{defi}

\begin{defi} Un paraboloide  hiperbólico es el conjunto de puntos en el espacio cuya distancia a un plano fijo  de manera proporcional al cuadrado de su distancia a dos líneas fijas en ese plano, pero con signos opuestos. Sus secciones paralelas a su eje de simetría son parábolas, mientras que las secciones perpendiculares son hipérbolas, dándole su característica forma de ``silla de montar".
\begin{equation}c(z-j)=\frac{(x-h)^2}{a^2}-\frac{(y-k)^2}{b^2}\end{equation}
\end{defi}
\begin{defi} Un hiperboloide de una hoja es el lugar geométrico de los  puntos en el espacio cuya diferencia de distancias a dos planos paralelos es constante.Se caracteriza por ser una superficie reglada, lo que significa que por cada punto de la superficie pasan dos rectas generatrices. Su ecuación general es \begin{equation}\frac{(x-h)^2}{a^2}+\frac{(y-k)^2}{b^2}-\frac{(z-j)^2}{c^2}=1.\end{equation}
\end{defi}

\begin{defi}Una hipérbole de dos hojas es el conjunto de puntos en el espacio cuya diferencia de distancias a dos planos paralelos es constante, pero con una separación mayor que en el caso de una hoja. Su ecuación canónica es 
\begin{equation}-\frac{(x-h)^2}{a^2}-\frac{(y-k)^2}{b^2}+\frac{(z-j)^2}{c^2}=1\end{equation} y esta compuesto por dos partes separadas, cada una con curvatura positiva en todos sus puntos.
\end{defi}

\subsection{Grupos de transformaciones}
Los grupos discretos $D_n$ (diedrales) y continuos $O(3)$ (ortogonales) son fundamentales para clasificar simetrías:

\begin{itemize}
    \item $D_n$: Grupo de simetrías de un $n$-ágono regular (rotaciones $\rho_{2\pi/n}$ y reflexiones).
    \item $O(3)$: Grupo de matrices $3\times 3$ que preservan distancias (rotaciones y reflexiones en $\mathbb{R}^3$).
\end{itemize}

Estos grupos actúan naturalmente sobre cónicas y cuadráticas. Por ejemplo, una esfera y un elipsoide tiene grupo de simetría $O(3)$, de hecho, 
toda superficie cuadrática no degenerada en $\mathbb{R}^3$ es isométrica a uno de estos casos canónicos ya definidos.

\section{Grupos de simetría en $\mathbb{R}^2$ y $\mathbb{R}^3$}

\begin{defi}
Una figura $\mathscr{F}$ es simétrica con respecto a un punto $O\in \mathbb{R}^2$ si para cada punto $P\in\mathscr{F} $ existe un punto $P^*\in \mathscr{F}$ tal que\\
i) $ P,O,P^*$ están en la misma línea (son colineales) \\
ii) $d(P,O)=d(P^*,O)$\\ 
Entonces, $O$ es centro de simetría de $\mathscr{F}$
\end{defi} 
\begin{teo}\label{teo-punto}
Si $P=(p_1,p_2) \in \mathscr{F}\implies P^*=(-p_1,-p_2) \in \mathscr{F}$ si y solo si $\mathscr{F}$ tiene simetría central en  $O=(0,0)$
        
\end{teo}

\begin{proof}[Prueba]

La recta que contiene a los  segmentos  $\overline{OP}$ y $\overline{OP^*}$ tiene la misma pendiente
$$m= \frac{p_2-0}{p_1-0}=\frac{p_2}{p_1}$$
y por lo tanto son colineales. Con eso se puede probar la primera condición. Para probar la segunda basta con verificar $$d(P,O)= \sqrt{p_1^ 2+p_2 ^2}=
\sqrt{(-p_1)^2+(-p_2)^2}=d(P^*,O)$$ 
\end{proof}

\begin{ejem}
El centro de simetría de una circunferencia de radio $r$ coincide con su centro natural. Supongamos que $P=(p_1,p_2)$ pertenece a la circunferencia con centro en el origen y radio $r$, entonces
$$p_1^2+p_2^2=r^2\implies (-p_1)^2+(-p_2)^2=r^2$$
Por teorema \ref{teo-punto} la circunferencia con centro en el origen tiene simetría central con respecto a su centro. Usando traslaciones, podemos inferir que el centro de simetría de la circunferencia  (\ref{circun})
es $(h,k)$
\end{ejem}
Los puntos simétricos de una figura están a la misma distancia del eje de simetría en direcciones opuestas de manera que la recta que une al punto y su simétrico es perpendicular al eje de simetría.
Veremos dos maneras de definir la simetría axial, primero, mediante estas propiedades y luego veremos más adelante una definición usando transformaciones la cual denominaremos reflexión.\\
\begin{defi}
    
Una figura $\mathscr{F}$ es simétrica con respecto a una recta ${l}$ si para cada punto $P\in\mathbb{R}^2$ existe un $P^*\in\mathbb{R}^2$ tal que\\
i) $\overline{PP^*}\bot{}\ l$\\
ii) $d(P,{l})=d(P^*,{l})$\\
Luego ${l}$ es un eje de simetría y $\mathscr{F}$ tiene simetría axial respecto a ${l}$.
\end{defi}
\begin{teo}\label{simeaxialr2}
Si $P= (p_1, p_2)\in \mathbb{R}^2\implies P^*=(-p_1, p_2)\in \mathbb{R}^2$ si y solo si, la figura es simétrica con respecto al eje Y del plano cartesiano. \end{teo}

\begin{proof}[Prueba]
    La pendiente de la recta que contiene al segmento  $\overline{PP^*}$ es \\$$m= \frac{p_2-p_2} {p_1+p_1}=0$$ lo que indica que esta sobre una recta horizontal y es perpendicular al eje $Y$ \\
$$d(P,{l})=\left|\frac{(p_1)+0(p_2)+0}{ \sqrt{1^2+0}}\right|
\\ =\left|\frac{p_1}{1}\right|= \left|\frac{-p_1}{1}\right|=\left|\frac{(-p_1)+0(p_2)+0}{\sqrt{1^2+0^2}}\right| =d(P^*,{l})$$
\end{proof}

\begin{ejem}
    Un triángulo isósceles tiene simetría axial con respecto a la mediatriz de la base o lado desigual.  Consideremos el triángulo con vértices $(p_1,p_2)$ $(-p_1,p_2)$  y $(0,p_3)$ con $p_3\ne p_2$. Según el teorema \ref{simeaxialr2} la mediatriz que coincide con el eje $Y$ es una recta de simetría axial para el triángulo, porque para todo punto $P\in\mathscr{F}$ implica que $P^*\in\mathscr{F}$.
    
\end{ejem}
\begin{defi}
La simetría rotacional es una transformación lineal que manda al vector $(p_1,p_2)$ asociado al punto $(x,y)$ otro vector $(p_1^*,p_2^*)$ asociado a un punto final que resultará de la rotación de un $\theta$ (en sentido antihorario respecto al eje $X$ positivo) y de la matriz \\
\begin{equation}
\begin{pmatrix}
cos(\theta) & -sen(\theta) \\
sen(\theta) & cos(\theta)
\end{pmatrix}
\end{equation}
\end{defi}
Por  ejemplo, si rotamos un cuadrado con respecto a su centro (intersección de sus diagonales) con ángulos de $90^0$, $180^0$, $270^0$ y $360^0$, obtendríamos el mismo cuadrado. Supongamos que $P_1=(1,1)$ es la esquina superior derecha del cuadrado, hacemos una rotación de  $90^0$ ese punto se transformaría en el punto de la esquina superior izquierda mediante la operación
\begin{equation*}
\begin{pmatrix}
cos(90^0) & -sen(90^0) \\
sen(90^0) & cos(90^0)
\end{pmatrix}\begin{pmatrix}
1  \\
1 
\end{pmatrix}=\begin{pmatrix}
-1  \\
1 
\end{pmatrix}
\end{equation*}
y entonces $P^*=(-1,1)$. Lo mismo aplica para todos los puntos de la figura obteniendo el mismo cuadrado pero rotado a $90^0$.

Esta matriz tiene determinante igual a $1$, al contrario de la reflexión cuyo determinante es $-1$.

\begin{defi}
La matriz de reflexión respecto al eje $X$, es la transformación \begin{equation*}
\begin{pmatrix}
1 & 0 \\
0 & -1
\end{pmatrix}\end{equation*}
y es una reflexión con respecto al eje $Y$ cuando 
\begin{equation*}
\begin{pmatrix}
-1 & 0 \\
0 & 1
\end{pmatrix}\end{equation*}

\end{defi}

A continuación, veremos como luciría lo anterior en $\mathbb{R}^3$. Se presentarán tres definiciones, junto con tres teoremas asociados.

\begin{defi}
Sea $\mathbb{S}$ una superficie y $O$ un punto en el espacio. Entonces $\mathbb{S}$ tiene simetría central con respecto a $O$ si y solo si, para cada $P\in\mathbb{S}$ existe un $P^*\in\mathbb{S}$ tal que\\
i) $O,P,P^*$ son colineales.\\
ii) $d(O,P)=d(O,P^*)$ \end{defi}
\begin{teo}\label{simecentralr3}
Si $P=(p_1,p_2,p_3) \in\mathbb{S}\implies P^*(-p_1,-p_2,-p_3)\in\mathbb{S}$ si y solo si $O=(0,0,0)$ es centro de simetría de $\mathbb{S}$.
\end{teo}
\begin{proof}[Prueba]

i) Sea la recta 
$L:(x,y,z)=(0,0,0)+t(p_1,p_2,p_3)$ Evaluando t en $0$, $1$, $-1$ se demuestra que $O,P,P^*$ pertenecen a $L$ y por lo tanto son colineales.

\begin{equation*}
ii)\ d(O,P)=\sqrt{(p_1-0)^2+(p_2-0)^2+(p_3-0)^2} \end{equation*}
\begin{equation*}
= \sqrt{(-p_1-0)^2+(-p_2-0)+(-p_3-0)^2}=d(O,P^*)\end{equation*}

\end{proof}
\begin{ejem}
    Los  paraboloides  no tienen centro de simetría. Sea $\mathbb{S}$ un paraboloide elíptico. Supongamos que $P=(p_1,p_2,p_3)\in \mathbb{S}$, por \ref{paraelip} se tiene  
$$c(p_3-j)=\frac{(p_1-h)^2}{a^2}+\frac{(p_2-k)^2}{b^2}$$
Realizamos una traslación al origen, de tal manera que tenemos una superficie $\mathbb{S^*}$ igual a
$$cp_3=\frac{p_1^2}{a^2}+\frac{p_2^2}{b^2}=\frac{(-p_1)^2}{a^2}+\frac{(-p_2)^2}{b^2}=-cp3$$
Lo cual es una contradicción. Luego, $P^*=(-p_1,-p_2,-p_3)\notin\mathbb{S^*}$.
Por teorema \ref{simecentralr3}, $\mathbb{S^*}$ no tiene que simetría central con respecto al $(0,0,0)$ y al trasladarnos implicaría que $(h,k,j)$ no es un punto de simetría central de $\mathbb{S}$. De hecho el exponente impar de una de las variables de las superficies garantiza la no existencia de puntos de simetría central.
\end{ejem}
\begin{defi}
Sea $\mathbb{S}$ una superficie y $l$ una recta en el espacio. Entonces $\mathbb{S}$ tiene simetría axial con respecto a $l$ si y solo si para cada punto $P\in \mathbb{S}$ existe $P^*\in\mathbb{S}$ tal que:\\

\noindent i) $\overline{PP}\perp l$\\
ii) $d(P,l)=d(P^*, l)$
\end{defi}
\begin{teo}\label{simeaxial3d} Si $P=(p_1,p_2,p_3) \in\mathbb{S}\implies P^*=(p_1,-p_2,-p_3)\in\mathbb{S}$ si y solo si  $\mathbb{S}$ es simétrico con respecto al eje $X$.
\end{teo}
\begin{proof}[Prueba]

i) La recta que contiene al segmento $\overline{PP}$ viene dada por $$l: (x,y,z)= ( p_1,p_2,p_3)+t(p_1-p_1,p_2-p_2,p_3-p_3).$$ Por otra parte,
el eje $X$ está dada por la recta $(x,y,z)= (0,0,0)+t(1,0,0)$. Por lo tanto, el producto punto de sus vectores 
directores es cero. Luego, son ortogonales.\\

ii) La distancia entre $P$ y la recta  y la recta y $P^*$ son las mismas. Llamemos $\vec{v}$ al vector director de la recta o eje $X$, 
$$d(P,l)=\frac{\left\|\overline{OP}\times \vec{v}\right\|}{\left\|\overline{v}\right\|}=\frac{\left\|det\left(\begin{array}{ccc}
\vec{i} & \vec{j} & \vec{k} \\
p_1 & p_2 & p_3 \\
1 & 0 & 0
\end{array}\right)\right\|}{\sqrt{1^2+0^2+0^2}}$$ $$=\|(0,p_3,-p_2)\|=\|(0,-p_3,p_2)\|=d(P^*,l)$$
\end{proof}

\begin{defi}
    Sea $\mathbb{S}$ una superficie y $\pi$ un plano. Entonces $\mathbb{S}$ es simétrica con respecto al plano si y solo si para cada $P\in\mathbb{S}$ existe un $P^*\in \mathbb{S}$ tal que\\
    i) $\overline{PP^*}$ es ortogonal a $\pi$\\
    ii) $d(P,\pi)=d(P^*,\pi)$
\end{defi}

\begin{teo}\label{simeplano3d}
Si $P=(p_1,p_2,p_3)\in \mathbb{S}\implies P^*=(-p_1,p_2,p_3)\in\mathbb{S}$ si solo si $\mathbb{S}$ es simétrica con respecto al plano $YZ$, o en otras palabras el plano en su ecuación general $x=0$.
\end{teo}
\begin{proof}[Prueba]
    
i) El vector director de la recta que contiene a los puntos $P$ y $P^*$ es $(2p_1,0,0)$ y el del plano es $(1,0,0)$, obteniendo cero al calcular el producto punto de ambos. Se concluye que $\overline{PP^*}$ es ortogonal a $\pi$.\\
ii) La norma del vector director del plano $x=0$ es $1$. 
$$d(P,x=0)=\left\|\frac{1(p_1)+0(p_2)+0(p_3)+0}{1}\right\|=|p_1|=|-p_1|=d(P^*,x=0)$$
\end{proof} 
\begin{ejem}
    La superficie cuadrática escrita en su forma general $$48x^2+32y^2-24z^2+96x-320y-960z-8994=0$$ tiene eje de simetría $$l:(x,y,z)=(-1,5,-20)+t(0,0,1)$$ y por otra parte $x=-1$, $y=5$  $z=-20$ son planos de simetría del hiperboloide de una hoja.\end{ejem}  
    Primero procedemos a escribir la superficie en su forma estándar o canónica factorizando y completando cuadrados.
    $$48(x^2+2x)+32(y^2-10y)-24(z^2+40z)=8994$$
    $$48(x+1)^2+32(y-5)^2-24(z+20)^2=8994+48(1)^2+32(5)^2-24(20)^2=192$$
    Dividimos todo en $192$ obteniendo
$$\frac{(x+1)^2}{4}+\frac{(y-5)^2}{6}-\frac{(z+20)^2}{8}=1$$
    usando una traslación al origen y los  teoremas \ref{simeaxial3d} y \ref{simeplano3d} podemos garantizar lo requerido puesto que los cuadrados en las variables implica que si colocamos algún negativo, el cuadrado lo anularía.\\
    
    En efecto, trasladamos al origen. Obtenemos
    \begin{equation}\label{hiperori}
    \frac{x^2}{4}+\frac{y^2}{6}-\frac{z^2}{8}=1\end{equation}
    Si suponemos que $P=(p_1,p_2,p_3)$ pertenece al hiperboloide trasladado (\ref{hiperori}) se verifica que $P^*=(-p_1,-p_2,p_3)$ también pertenece y, por consiguiente, el eje $Z$ es un eje de simetría. De igual manera  $P^*=(-p_1,p_2,p_3)$, $P^*=(p_1,-p_2,p_3)$ y $P^*=(p_1,p_2,-p_3)$  implicaría que los planos $x=0$, $y=0$, $z=0$ son planos de simetría del hiperboloide trasladado. Regresandonos de la traslación se obtiene lo deseado.
\section{Modelado matemático de estructuras arquitectónicas}

\subsection{Panteón de Agripa: Simetría $D_{8}$}
El Panteón de Agripa, reconstruido por Adriano (118-125\,d.C.), es paradigmático por su simetría estructural y carga simbólica. Su cúpula, de 43.3\,m de diámetro con casetones concéntricos y un óculo central, no solo optimiza la distribución de tensiones mediante una geometría octogonal implícita, sino que también proyecta patrones lumínicos invariantes bajo rotaciones, encarnando el ideal romano de orden cósmico \cite{HNG, Caballero2025}. Esta disposición, modelable mediante el grupo diédrico \(D_{8}\), será la base para nuestro análisis algebraico.
La cúpula del Panteón  idealiza un octágono regular en planta, cuyo grupo de simetría  explica la distribución equitativa de tensiones estructurales.\\
\begin{figure}
    \centering
    \includegraphics[width=0.75\linewidth]{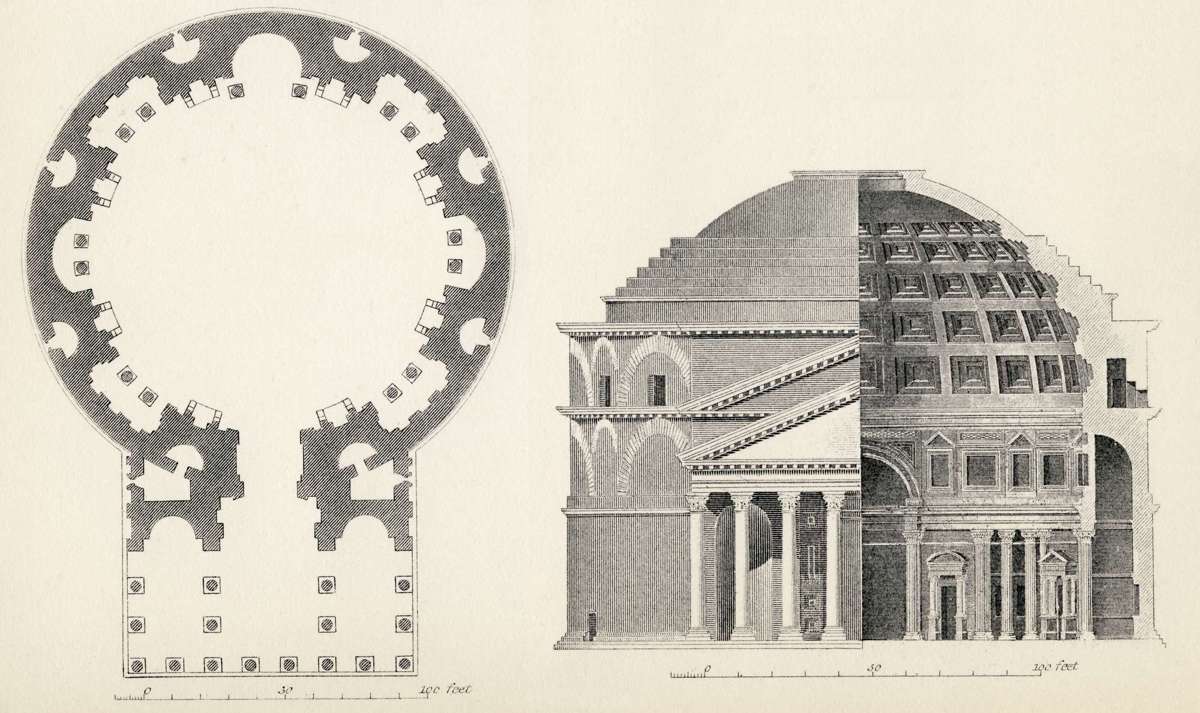}
    \caption{Panteón de Agripa.\cite{HNG} }
    \label{fig:placeholder}
\end{figure}

Matricialmente, las generadoras son:
\begin{equation}
\rho_{\pi/4} = \begin{pmatrix}
\cos\frac{\pi}{4} & -\sin\frac{\pi}{4} & 0 \\
\sin\frac{\pi}{4} & \cos\frac{\pi}{4} & 0 \\
0 & 0 & 1
\end{pmatrix}, \quad
\sigma_h = \begin{pmatrix}
1 & 0 & 0 \\
0 & 1 & 0 \\
0 & 0 & -1
\end{pmatrix}
\end{equation}
Para cuantificar la ventaja estructural que confiere la simetría, se implementó un modelo numérico simplificado mediante el Método de Elementos Finitos (MEF) explicado en \cite{logan}. Se comparó el comportamiento de un pórtico planar bajo carga gravitacional y lateral, en dos configuraciones: una simétrica (\(D_2\)) y otra asimétrica.

Cada barra $e=(i,j)$ une los nudos $i$ y $j$ con longitud $L_e$ y dirección $(c,s)=(\cos\theta,\sin\theta)$.
El elemento es axial, por lo que su matriz de rigidez en coordenadas globales es
\begin{equation}
k_e \;=\; \frac{EA}{L_e}
\begin{bmatrix}
 c^2 & cs & -c^2 & -cs \\
 cs  & s^2 & -cs  & -s^2 \\
 -c^2& -cs & c^2  & cs \\
 -cs & -s^2& cs   & s^2
\end{bmatrix},
\qquad
(c,s)=\left(\frac{\Delta x}{L_e},\frac{\Delta y}{L_e}\right),
\label{eq:ke}
\end{equation}
donde $E$ es el módulo de elasticidad y $A$ el área de la sección transversal. La rigidez global $K$ se obtiene por ensamblaje estándar de (\ref{eq:ke}) sobre los grados de libertad de cada nudo. El sistema lineal a resolver es
\begin{equation}
K\,u \;=\; F,
\label{eq:ku_f}
\end{equation}
con $u$ el vector de desplazamientos nodales y $F$ el vector de cargas nodales. Las reacciones aparecen al recuperar los esfuerzos en apoyos. Para cada barra, el esfuerzo axial (convención: $+$ tracción, $-$ compresión) se calcula como
\begin{equation}
N_e \;=\; \frac{EA}{L_e}\,[-c,\,-s,\,c,\,s]\;\,\cdot u_e,
\label{eq:axial}
\end{equation}
siendo $u_e$ el vector de desplazamientos de los dos nudos extremos del elemento, reordenado en el sistema global.\\

Se consideran cuatro nudos $(A,B,C,D)$ formando un marco rectangular en el caso simétrico $D_2$, y una traslación del pilar izquierdo en el caso asimétrico:
\begin{align*}
\textbf{Caso simétrico ($D_2$):}\quad
&A(-1,0),\; B(-1,2),\; C(1,2),\; D(1,0).\\[2pt]
\textbf{Caso asimétrico:}\quad
&A(-0.5,0),\; B(-0.5,2),\; C(1,2),\; D(1,0).
\end{align*}

Se usan las barras perimetrales y dos diagonales (arriostramiento en X) para evitar mecanismos:
\[
E \;=\; \{\,AB,\;BC,\;CD,\;DA,\;AC,\;BD\,\}.
\]

Se considera los apoyos traslacionales en la base como
\[
u_A = u_D = (0,0).
\]
 y las cargas nodales
\[
F_y(A)=F_y(B)=F_y(C)=F_y(D)=-500,\qquad F_x(B)=+1000,
\]
que representan, respectivamente, peso propio simplificado y viento en el nudo B.

Supongamos además,
\[
E=210\;\text{GPa},\qquad A=1\;\text{cm}^2=1\times10^{-4}\;\text{m}^2.
\]

En la Tabla~\ref{tab:geom-sim} y \ref{tab:geom-asim} se listan $L_e$ y $(c,s)$ para sustituir directamente en (\ref{eq:ke}).
\begin{table}[H]
\centering
\caption{Caso simétrico $D_2$.}
\label{tab:geom-sim}
\begin{tabular}{c c c c}
Elemento & $L_e$ (m) & $c=\cos\theta$ & $s=\sin\theta$\\\hline
$AB$ & $2.0000$ & $0.0000$ & $+1.0000$\\
$BC$ & $2.0000$ & $+1.0000$ & $0.0000$\\
$CD$ & $2.0000$ & $0.0000$ & $-1.0000$\\
$DA$ & $2.0000$ & $-1.0000$ & $0.0000$\\
$AC$ & $2.8284$ & $+0.7071$ & $+0.7071$\\
$BD$ & $2.8284$ & $+0.7071$ & $-0.7071$\\
\end{tabular}
\end{table}

\begin{table}[H]
\centering
\caption{Caso asimétrico.}
\label{tab:geom-asim}
\begin{tabular}{c c c c}
Elemento & $L_e$ (m) & $c=\cos\theta$ & $s=\sin\theta$\\\hline
$AB$ & $2.0000$ & $0.0000$ & $+1.0000$\\
$BC$ & $1.5000$ & $+1.0000$ & $0.0000$\\
$CD$ & $2.0000$ & $0.0000$ & $-1.0000$\\
$DA$ & $1.5000$ & $-1.0000$ & $0.0000$\\
$AC$ & $2.5000$ & $+0.6000$ & $+0.8000$\\
$BD$ & $2.5000$ & $+0.6000$ & $-0.8000$\\
\end{tabular}
\end{table}
Se ensambló $K$ para el caso simétrico a partir de (\ref{eq:ke}), obteniendose
$$
\begin{pmatrix}
14212346 & 3712346 & 0& 0& 3712346& 3712346 & -10500000& 0 
\\
3712346 & 14212346 & 0 & -10500000& 3712346& 3712346 & 0& 0
\\
0 & 0 & 14212346& -3712346& -10500000& 0 & -3712346& 0
\\
0 & -10500000 & -3712346& 14212346& 0& 0 & 3712346& -3712346 
\\
3712346 & -3712346 & -10500000 & 0& 14212346& 3712346 & 0& 0
\\
3712346 & 3712346 & 0& 0& 3712346& 14212346 & 0& -10500000
\\
-10500000 & 0 & -3712346 & 3712346& 0& 0 & 14212000& -3712346
\\
0 & 0 & 3712346& -3712346& 0& -10500000 & -3712346 & 14212346
\end{pmatrix}$$
se impusieron las restricciones en $A$ y $D$ y se resolvió (\ref{eq:ku_f}) para los cuatro grados de libertad libres $\{u_{Bx},u_{By},u_{Cx},u_{Cy}\}$. Luego,
$$\begin{pmatrix}
 14212346& -3712346& -10500000& 0 \\
 -3712346& 14212346& 0& 0\\
-10500000 & 0& 14212346& 3712346\\
0& 0& 3712346& 14212346
\end{pmatrix}\begin{pmatrix}
 u_{Bx} \\
 u_{By}\\
u_{Cx} \\
u_{Cy}
\end{pmatrix}=\begin{pmatrix}
 1000 \\
 -500\\
0 \\
-500
\end{pmatrix}$$

Los desplazamientos se reportan en milímetros (mm) y los esfuerzos axiales de barra a partir de (\ref{eq:axial}) en Newtons (N). Por lo tanto, $$\begin{pmatrix}
 u_{Bx} \\
 u_{By}\\
u_{Cx} \\
u_{Cy}
\end{pmatrix}=\begin{pmatrix}
 0.1979 \\
 0.0165\\
 0.1667\\
-0.0787
\end{pmatrix}$$

Los desplazamientos nodales (mm) son
\[\textbf{D}_2: \hspace{1cm}  B:(0.1979,\;0.0165), \hspace{1cm}  C:(0.1667,\;-0.0787)
\]
\\
De igual manera, se puede verificar
\[ \textbf{Asim.:} \hspace{1cm} B:(0.2624,\;0.0352), \hspace{1cm} C:(0.2376,\;-0.0919).
\]

Obteniendo una norma euclídea de desplazamientos como indicador de deformabilidad global
\[
\|\mathbf{u}\|_2 \;=\; \mathbf{0.271}\ \text{mm} \quad (\text{D}_2), 
\qquad
\|\mathbf{u}\|_2 \;=\; \mathbf{0.367}\ \text{mm} \quad (\text{asim.}).
\]

Los esfuerzos axiales por barra son
\[
\begin{array}{l rrrrrr}
 & AB & BC & CD & DA & AC & BD \\\hline
\textbf{D}_2 & +173 & -328 & -826& \phantom{-}0 & +462 & \mathbf{-952}\\
\textbf{Asim.} & +370 & -347 & \mathbf{-964} & \phantom{-}0 & +580 & \mathbf{-1086}
\end{array}
\]
En negrita se señalan los picos por caso; se observa que la asimetría incrementa el máximo absoluto de esfuerzo, por ejemplo, en la barra $BD$ ($\approx 14\%$ más) y en la deformación global ($\approx 35\%$ más) bajo las mismas cargas y propiedades.

La configuración $D_2$ redistribuye las tensiones, los caminos de carga simétricos reducen los grados de libertad flexibles del sistema ensamblado y limitan la aparición de concentraciones de esfuerzo. La rotura de simetría desplazamiento del pilar izquierdo induce trayectorias excéntricas de carga: el arriostramiento diagonal derecho recibe mayor demanda y la norma de desplazamientos aumenta. Esta evidencia numérica respalda que la simetría no es sólo una propiedad geométrica elegante, sino un principio estructural que tiende a optimizar rigidez y distribución de tensiones.

El salto conceptual a $D_8$ consiste en replicar el módulo $D_2$ por rotación de $45^\circ$ en planta. En un anillo de ocho sectores, el ensamblaje de rigideces exhibe bloques rotacionalmente equivalentes y condiciones periódicas que suavizan los desplazamientos relativos entre sectores y aplanan picos de esfuerzo. Esta modularidad por simetría es consistente con estrategias de modelado numérico para edificaciones históricas, donde el uso de periodicidad y subestructuras reduce incertidumbres locales y mejora la estabilidad global \cite{sala2024}. En términos del grupo, tanto las rotaciones como las reflexiones de $D_8$ restringen las formas admisibles de $u$ en (\ref{eq:ku_f}).
\begin{figure}
    \centering
    \includegraphics[width=0.75\linewidth]{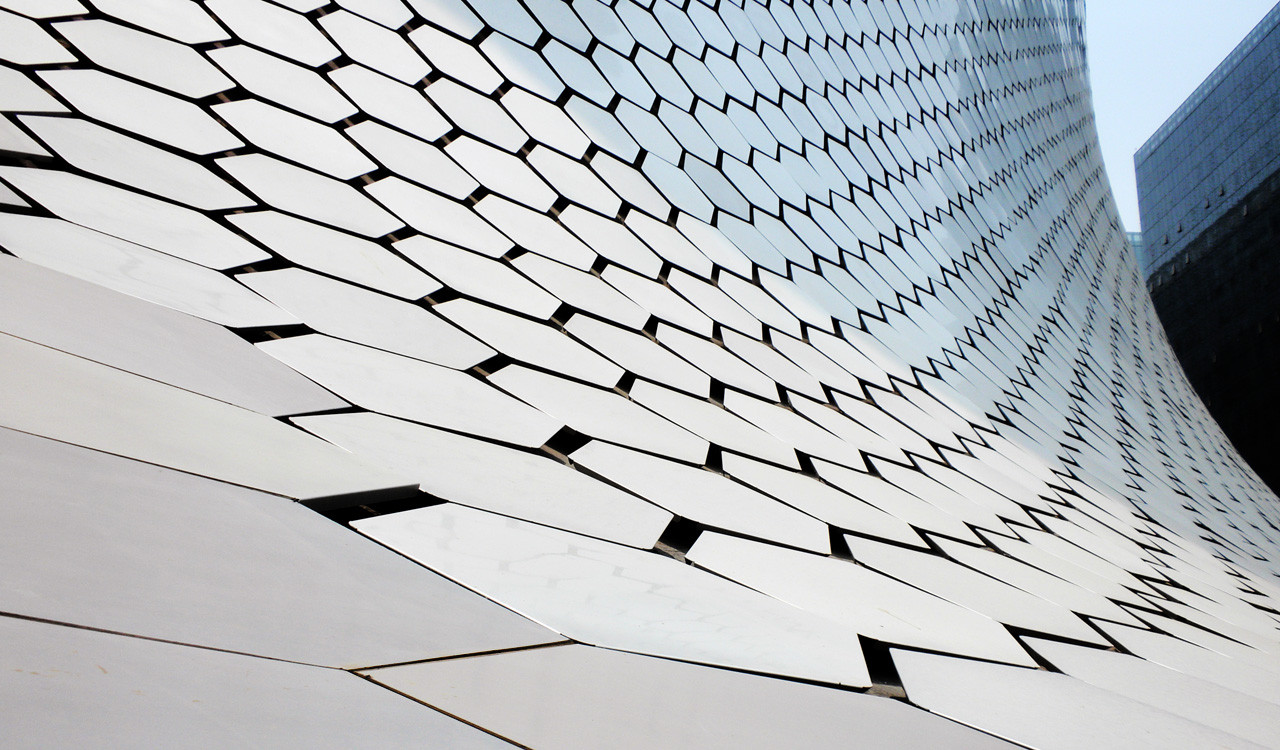}
    \caption{Fachada del museo Soumaya, México \cite{ref11}}
    \label{fig:placeholder}
\end{figure}
Por otra parte, sería interesante profundizar sobre cómo las características simétricas del Museo Soumaya de México no sólo la hace hermosa a la vista sino que también da otro plus. Su hexarotación en fachada \cite{ref11} se determina por la matriz

    \[
    \rho_{2\pi/6} \cdot \mathbf{v} = \begin{pmatrix}
    \cos\theta & -\sin\theta & 0 \\
    \sin\theta & \cos\theta & 0 \\
    0 & 0 & 1
    \end{pmatrix} \mathbf{v}, \quad \theta = \frac{2\pi}{6}
    \]
     y  optimiza la conservación y durabilidad de todo el edificio \cite{ref11}.

En conclusión, este artículo indica que el estudio algebraico de las simetrías, mediante grupos de transformaciones como $D_n$ y
$O(3)$, trasciende el ámbito teórico para constituir una herramienta poderosa en el análisis geométrico y estructural. La clasificación de cónicas y superficies cuadráticas a través de sus invariantes sienta las bases para un marco analítico riguroso. La aplicación de este marco al Panteón de Agripa y al Museo Soumaya revela que la simetría axial y rotacional no es un mero recurso estético, sino un principio fundamental de eficiencia estructural y funcional. El modelo de elementos finitos confirma cuantitativamente que la configuración 
$D_8$
  del Panteón optimiza la distribución de tensiones y limita las deformaciones, mientras que la simetría hexagonal del Soumaya sugiere una optimización lumínica. Así, se establece un puente indisociable entre la abstracción del álgebra, la precisión de la geometría y la innovación en el diseño arquitectónico, abriendo puertas a futuras investigaciones interdisciplinarias.

\label{end-art}
\end{document}